\documentclass[12pt, a4paper]{amsart}

\usepackage{amsthm,amssymb,anysize, hyperref, txfonts, dsfont, enumerate, color,booktabs}

\marginsize{2cm}{2cm}{2cm}{2cm}

\theoremstyle{plain}
\newtheorem{theorem}{Theorem}[section]
\newtheorem{lemma}[theorem]{Lemma}

\newcommand\PG{\mathrm{PG}}
\newcommand\F{\mathbb{F}}
\newcommand\Aut{\mathrm{Aut}}
\newcommand{\gaussian}[2]{\genfrac{[}{]}{0pt}{}{#1}{#2}}
\newcommand\ones{\mathds{1}}

\newcommand\charvec{\ones}

\renewcommand\ge{\geqslant}
\renewcommand\le{\leqslant}

\title[On the $430$-cap of $\PG(6,4)$]{On the $430$-cap of $\PG(6,4)$ having two\\ intersection sizes with respect to hyperplanes}

\author{John Bamberg}

\email{John.Bamberg@uwa.edu.au}

\address{ %
Centre for the Mathematics of Symmetry and Computation\\
Department of Mathematics and Statistics\\
The University of Western Australia\\
Crawley, W.A. 6009, Australia.}


\begin{document}

\begin{abstract}
Let $\mathcal{C}$ be a 430-cap of $\PG(6,4)$ having two intersection sizes with respect to hyperplanes.
We show that no hyperplane of $\PG(6,4)$ intersects $\mathcal{C}$ in a Hill 78-cap.
So if it can be shown that the Hill 78-cap of $\PG(5,4)$ is projectively unique, then
such a 430-cap does not exist, or equivalently, a two-weight $[430,7]_{\mathbb{F}_4}$ linear code with
dual weight at least 4, does not exist.
\end{abstract}

\maketitle

\section{Introduction}

Uniformly packed codes generalise perfect codes, and the 1-error correcting examples have connections to strongly regular graphs, partial quadrangles,
and two-character sets in finite projective spaces (see \cite{Calderbank:1986wc}).
An $e$-error correcting code $C$ is \emph{uniformly packed} if spheres of radius $e+1$ about codewords cover the whole space, and vectors at distance $e$ from the $C$ are in $\lambda+1$ spheres while vectors at distance $e+1$ from the code are in $\mu$ spheres. For the case that $e=1$, a code $C$ is $1$-error correcting
if and only if the dual code $C^\perp$ has two nonzero weights. If $C$ has minimum distance at least 3, then $C^\perp$ is \emph{projective},
and gives rise to a \emph{two-character set} of projective space: a set of points $S$ such that there are only two values for the possible intersection size
of a hyperplane with $S$. 

There is also a connection with finite \emph{partial quadrangles}.
Partial quadrangles were introduced by Cameron \cite{Cameron:1975vc}
as a (finite) geometry of points and lines such that every two points are on at most one
line, there are $s+1$ points on
a line, every point is on $t+1$ lines and satisfying the following two important
properties: (i) for every point $P$ and every line $\ell$ not incident with $P$,
 there is at most one point on $\ell$ collinear with $P$;
(ii) there is a constant $\mu$ such that for every pair of non-collinear points $(X,Y)$ there are precisely $\mu$ points collinear with $X$ and $Y$.
The collinearity graph of a partial quadrangle is strongly
regular. The only known partial quadrangles, which are not generalised quadrangles, are
triangle-free strongly regular graphs, arise from removing points from a generalised quadrangle of order $(s,s^2)$,
or one of three \emph{exceptional examples} arising from linear
  representation of one of the Coxeter $11$-cap of $\PG(4,3)$, the Hill
  $56$-cap of $\PG(5,3)$ or the Hill $78$-cap of $\PG(5,4)$.

A \textit{$k$-cap} of a projective space $\PG(n,q)$ is a set of $k$
points with no three collinear.  Calderbank \cite{Calderbank:1982ty} proved using
number-theoretic arguments that if a partial quadrangle is a linear
representation of a $k$-cap then $q\ge 5$ or it is isomorphic to the linear representation of
one of the following: (i) an ovoid of $\PG(3,q)$; (ii) the Coxeter $11$-cap of $\PG(4,3)$;
(iii) the Hill $56$-cap of $\PG(5,3)$; (iv) a $78$-cap of $\PG(5,4)$; (v) a $430$-cap of $\PG(6,4)$.
Tzanakis and Wolfskill \cite{Tzanakis:1986wx} proved that
if $q\ge 5$, then the first case applies. It is still not known if case (v) occurs; that is, whether there is a $430$-cap of $\PG(6,4)$ 
such that every hyperplane intersects in 78 or 110 elements. If a hyperplane intersects in 78 elements, then it is a two-character 78-cap
of $\PG(5,4)$ (see Lemma \ref{twocharacter}). This leaves two open problems: 
\begin{enumerate}
\item Does there exist a two-character 78-cap of $\PG(5,4)$ projectively inequivalent to Hill's cap?
\item Does there exist a two-character 430-cap of $\PG(6,4)$?
\end{enumerate}
These problems are of interest to finite geometry and coding theory alike, and have been
open for over 40 years, since at least \cite{Calderbank:1982ty}. We show in this note that
a negative solution to the first problem implies a negative solution to the second problem.



\begin{theorem}\label{main}
Let $\mathcal{C}$ be a $430$-cap of $\PG(6,4)$ having two intersection sizes with respect to hyperplanes.
Then no hyperplane of $\PG(6,4)$ intersects $\mathcal{C}$ in a cap projectively equivalent to the Hill $78$-cap.
\end{theorem}

The basic argument proceeds as follows.
Suppose $H$ is the Hill $78$-cap of $\PG(5,4)$ and 
embed $\PG(5,4)$ as a hyperplane $\Pi$ of $\PG(6,4)$.
Let $\mathcal{Q}$ be the partial quadrangle arising from linear representation of $H$,
and let $\Gamma$ be its collinearity graph. Then $\Gamma$ is a strongly regular graph with
parameters $(4096,234,2,14)$.
Now the affine points are the points of $\mathcal{Q}$, and the affine lines meeting $\Pi$ in a point of $H$
are the lines of $\mathcal{Q}$.
 Let $\mathcal{C}$ be a $430$-cap of $\PG(6,4)$ containing $H$. 
 So the affine points $\bar{\mathcal{C}}:=\mathcal{C}\backslash \Pi$ of $\mathcal{C}$ form a set of points of size $352$ of $\mathcal{Q}$
 such that every line of $\mathcal{Q}$ intersects it in at most one point.
Moreover, $\bar{\mathcal{C}}$ forms a \emph{Delsarte coclique} for $\Gamma$; a coclique that has size attaining the Delsarte/Hoffman bound.
We will show that $\Gamma$ does not have a Delsarte coclique, which then shows that the Hill $78$-cap does
not extend to a $430$-cap of $\PG(6,4)$. To do this, we take the Schurian scheme for the automorphism group of $\Gamma$, which is a 9-class fission scheme for the natural 2-class scheme arising from $\Gamma$. We then use another 2-class fusion of this Schurian scheme to yield information on the inner distribution of a
putative Delsarte coclique.

\section{Some background}

Let $\Omega$ be a set, and let $A_0, A_1, \ldots, A_d$ be symmetric  $\{0,1\}$-matrices with rows and columns indexed by $\Omega$. Then $\mathcal{A}=(\Omega, \{A_0, A_1, \ldots, A_d\})$ is a $d$-class \emph{association scheme} if the following conditions hold:
\begin{enumerate}
    \item $A_0$ is the identity matrix $I$,
    \item $\sum_{i=0}^d A_i$ is the matrix with every entry equal to $1$,
    \item There exist constants $p_{ij}^k$ depending only on $i$, $j$, and $k$, such that 
    $A_i A_j = \sum_{k=0}^d p_{ij}^k A_k$.
\end{enumerate}
The matrices $A_0$, $A_1$, $\ldots$, $A_d$ are the \emph{adjacency matrices} of $\mathcal{A}$, and indeed, 
each $A_i$ is the adjacency matrix
of an undirected graph. A strongly regular graph $\Delta$ is essentially equivalent to a $2$-class association scheme, where $A_1$ and $A_2$
are the adjacency matrices for $\Delta$ and its complement.

It is well known that $\mathbb{R}^\Omega$ decomposes into $d+1$ simultaneous eigenspaces for the adjacency matrices of $\mathcal{A}$. Moreover, there are projection matrices $E_0, E_1, \ldots, E_d$ (the \emph{minimal idempotents}) onto each of these eigenspaces, such that 
\[
E_i = \sum_{j=0}^d Q_{ji} A_j,
\]
where $Q$ is called the \emph{matrix of dual eigenvalues}.
If $C$ is a subset of $\Omega$, then its \emph{inner distribution} is the vector $a = (a_0,a_1,\ldots,a_d)$ defined by
\[
a_i= \frac{1}{|C|}\charvec_C A_i \charvec_C^\top.
\]
where we use $\charvec_C$ to denote the characteristic function of $C$ in $\Omega$.
If $Q$ is the matrix of dual eigenvalues of $\mathcal{A}$, then
\[
(a Q)_j = \frac{|\Omega|}{|C|} \charvec_C E_j \charvec_C^\top
\]
for all $j\ge 0$. The vector $aQ$ is sometimes known as the \emph{MacWilliams transform} of $C$,
and it follows from the fact that the $E_j$ are positive semidefinite, that 
 that each entry of $aQ$ is non-negative.
 
The \emph{outer distribution} $B$ of $S$ is the $|\Omega|\times d$ matrix, with rows indexed by $\Omega$ and columns indexed by the $R_i$, defined by
\[
B_{x,i}=|\{ y\in S : (x,y)\in R_i\}|=\charvec_{\{x\}}A_i\charvec_S^\top.
\]

The \emph{dual degree set} of $C$ is the set of nonzero indices $j$
for which the $j$-th coordinate of its MacWilliams transform is nonzero.
Two subsets of $\Omega$ are \emph{design-orthogonal} if
their dual degree sets are disjoint. In this case, we have the following elementary result, due at least to Roos. 

\begin{theorem}[{\cite[Corollary 3.3]{roos}}]\label{roos}
If $S,T\subset \Omega$ are design-orthogonal, then 
$|S\cap T|=\dfrac{|S|\cdot |T|}{|\Omega|}$.
\end{theorem}

A transitive group $G$ acting on $\Omega$ is \emph{generously transitive} if for any distinct pair $(\alpha,\beta)$
of elements of $\Omega$, there is some $g\in G$ such that $\alpha^g=\beta$ and $\beta^g=\alpha$.
If a finite group $G$ acts generously transitively on a set $\Omega$, then the orbits of $G$ on unordered pairs 
of $\Omega$ give rise to an association scheme.
We refer the reader to \cite{Godsil:2016vd} or \cite{Lint:2001ur} for more information on association schemes.

\section{A Schurian scheme and some interesting subsets}

The following \emph{cyclotomic} construction of $\Gamma$ can be found as 
\cite[Example FE3]{Calderbank:1986wc}.
Let $z$ be a primitive element of $\F_{4^6}$.
Let $O$ be $\langle z^{35}\rangle \cup \langle z^{35}\rangle z^7$.
Then $\Gamma$ is isomorphic to the Cayley graph $\mathrm{Cay}(V,O)$ where $V$ is the additive group of $\F_{4^6}$, and it is
strongly regular graph with parameters $(4096,234,2,14)$.
Note that we can also view $O$ as the set of underlying vectors of the Hill 78-cap (n.b., $234=3\times 78$), represented as elements
of $\F_{4^6}$. 

Some of the details below were aided by computer, and in particular, the \textsf{AssociationScheme} package \cite{AssociationSchemes} in \textsf{GAP}.
The automorphisms of $\Gamma$ are generated by the translations (addition by elements of $V$), multiplication by $z^{35}$, and the map $\rho:x\mapsto z^{42} x^4$ (of order 6). So $\Aut(\Gamma)$ is isomorphic to $C_2^{12}:(C_{117}:C_6)$.
Indeed, the stabiliser of $0$ is generated by $z^{35}$ and $\rho$, and 
these automorphisms act on $O$. Moreover, $\Aut(\Gamma)$ acts generously transitively
on the points of $\Gamma$. 

Let $v:=4096$, the number of vertices of $\Gamma$.
Take the Schurian association scheme $\mathcal{A}$ for $\Aut(\Gamma)$, which is a fission scheme
for the 2-class association scheme $\mathcal{G}$ associated to the original strongly regular graph $\Gamma$.
Then the valencies of $\mathcal{A}$ are (in order) 1, 117, 234, 234, 351, 351, 702, 702, 702, 702
with $R_2$ being the adjacency relation for $\Gamma$ (and the $R_i$ are indexed with $i\in\{0,\ldots,9\}$). 
In fact, $\mathcal{A}$ is a translation scheme and it is formally dual.
The matrix $P$ of eigenvalues, and the matrix $Q$ of dual eigenvalues, for $\mathcal{A}$ are:
\[
P=Q= 
\begin{bmatrix}
1&  117&  234&  234&  351&  351&  702&  702&  702&  702 \\
1&  -27&   10&   10&   15&   63&   30&   30&  -66&  -66 \\
1&    5&  -22&   10&   15&  -33&   30&   30&   30&  -66 \\
1&    5&   10&  -22&   15&  -33&   30&   30&  -66&   30 \\
1&    5&   10&   10&   47&   -1&  -34&  -34&   -2&   -2 \\
1&   21&  -22&  -22&   -1&   31&   -2&   -2&   -2&   -2 \\
1&    5&   10&   10&  -17&   -1&   30&  -34&   -2&   -2 \\
1&    5&   10&   10&  -17&   -1&  -34&   30&   -2&   -2 \\
1&  -11&   10&  -22&   -1&   -1&   -2&   -2&   30&   -2 \\
1&  -11&  -22&   10&   -1&   -1&   -2&   -2&   -2&   30 
\end{bmatrix}
  \]
We note that there is an involution normalising the association scheme. It is induced by the following semilinear map of order 12:
\[
\tau:x\mapsto z^{14} x^2
\]
and it interchanges relations of the association scheme: 
$R_2^\tau=R_3$, $R_6^\tau=R_7$, and $R_8^\tau=R_9$. 

There are some unions of relations in $\mathcal{A}$ that yield interesting graphs. 
\begin{enumerate}[(i)]
\item $R_2$ is the original strongly regular graph $\Gamma$. Moreover, the non-principal minimal idempotents for $\Gamma$
are $\sum_{j\in\{1, 3, 4, 6,7, 8\}}E_j$ and $\sum_{j\in\{2,5,9\}}E_j$, 
where $E_j$ is the $j$-th minimal idempotent for $\mathcal{A}$.

\item $R_2\cup R_7\cup R_8$ yields a strongly
regular Cayley graph $\mathcal{F}$ that will feature in our proof of Theorem \ref{main}. 
The elements of the subfield $\F_{4^3}$ form a maximal clique of $\mathcal{F}$.
\end{enumerate}

Below we list some interesting (Delsarte) designs for $\mathcal{A}$. We denote by $V_j$ the $j$-th eigenspace for $\mathcal{A}$,
for which the minimal idempotent $E_j$ projects to.

%

\subsection{Example 1: a subfield design}\label{subfield}

Consider the elements $U$ of $V$ that lie in the subfield $\F_{4^3}$. 
It turns out that the inner distribution of $U$ is $(1,0, 9, 0,  0, 0, 0, 27, 27, 0)$,
and so its MacWilliams transform is $(64, 0, 576, 0, 0, 0, 0, 1728, 1728, 0)$.
Therefore, 
$\charvec_U\in V_0\perp V_2\perp V_7\perp V_8$.

\subsection{Example 2: a Delsarte coclique}\label{Dcoclique}

The complement of $\Gamma$ is $k$-regular with $k=3861$, and it has least eigenvalue $\tau:=-11$.
The \emph{Delsarte bound} for the size of a coclique of $\Gamma$ is then $1-k/\tau=352$.
Suppose there exists a coclique $S$ of $\Gamma$ of size 352. Then $S$ is a \emph{Delsarte coclique}, and 
$\charvec_SE=0$ where $E:=\sum_{j\in\{1, 3, 4, 6,7, 8\}}E_j$ (see \cite[Corollary 3.7.2]{Godsil:2016vd}). (Note that $P_{i2}=-22$ for $i\in\{2,5,9\}$.)  
Recall that if $E:=\sum_{j\in\{1, 3, 4, 6,7, 8\}}E_j$, then $\charvec_SE=0$.
Consider $E_j$ where $j\in\{1, 3, 4, 6,7, 8\}$. Then $E_j=EE_j$ and so $\charvec_SE_j=\charvec_SEE_j=0$.
So we have $(aQ)_j=0$ for $j\in\{1, 3, 4, 6,7, 8\}$, or in other words,
\[
\charvec_S\in V_0\perp V_2\perp V_5\perp V_9.
\]
If we apply the involution $\tau$, we find that $\charvec_{S^\tau}\in V_0\perp V_3\perp V_5\perp V_8$.

Now consider the vector $v_P:=22\charvec_{\{P\}}+\charvec_{P^\perp}$ where $P$ is a point of $\mathcal{Q}$, and $P^\perp$
is the set of points adjacent to $P$. Notice that $v_P=\charvec_{\{P\}}(A_2+22I)$ and so $v_PE_j=0$ for $j\in\{2,5,9\}$.
In particular, $v_P$ is design-orthogonal\footnote{Note that design-orthogonality extends to weighted subsets in a straight-forward way.} to $S$
and so $\charvec_S\cdot v_P=22$. It follows that $|P^\perp\cap S|$ is equal to 22 when $P\notin S$, but equal to 0 when $P\in S$.

\section{Proof of Theorem \ref{main}}


\begin{proof}  
Let $S$ be a Delsarte coclique for $\Gamma$ and let $B$ be the outer distribution of $S$. 
By a theorem of Delsarte \cite[Theorem 3.1]{delsarte}, for all vertices $x$ of $\Gamma$, and for all $j\in\{1, 3, 4, 6, 7, 8\}$,
\begin{equation}
\sum_{i=0}^d \frac{P_{ji}}{P_{0i}} B_{x,i}=0.\label{eq:outer}
\end{equation}
Fix an element $x\notin S$. Recall that there are 22 elements of $S$ adjacent to $x$ (see \ref{Dcoclique} above), and so we can write
$B_{x,i}=(0, y_1, 22, y_2, y_3, y_4, y_5, y_6, y_7,   330 - y_1 - y_2 - y_3 - y_4 - y_5 - y_6 - y_7)$ for some $y_i$.
Then, we have the following equations (arising from \eqref{eq:outer}):
\begin{center}
\begin{tabular}{ll}
\toprule
$j$ & Equation\\
\midrule
1 & $220 + y_1 - (y_3 + y_4 + 2 y_5 + y_6 + y_7) = 0$\\
3 & $y_3 + y_5 + y_8 = 110$\\
4 & $y_1 + y_3 + 3 y_4 - y_6 - y_7=0$\\
6 & $y_1 + y_3 + y_6 - y_4 - y_7=0$\\
7 & $y_1 + y_3 + y_7 - y_4 - y_6=0$\\
8 & $2 (y_1 + y_3) - y_8=0$\\
\bottomrule
\end{tabular}
\end{center}

 
These equations reduce:
$2y_4=y_6=y_7 = y_8$,\quad $y_3 =110 - y_5 - y_8$, \quad $y_1 = y_5 + \tfrac{3}{2} y_8-110$.
So
\[
B_{x,i}=(0,  y_5 + \tfrac{3}{2} y_8-110, 22, 110 - y_5 - y_8, \tfrac{1}{2}y_8, y_5, y_8, y_8, y_8,  330 - y_5 - 4 y_8).
\]

Let $a$ be the inner distribution of $S$.
Now $A_2$ is the adjacency matrix of $\Gamma$, and so $\charvec_{S}A_2\charvec_{S}^\top=0$. Hence we can write the inner distribution of $S$ as 
$a=(1,x_1,0,x_3,x_4,x_5,x_6,x_7,x_8,351-x_1-x_3-x_4-x_5-x_6-x_8)$, where the $x_i$ are indeterminate.
Now multiply by $Q$ to yield the MacWilliams transform of $S$:
\begin{align*}
aQ=32(&11, \quad
\tfrac{1}{2}(- x_1 + x_3 + x_4 + 2 x_5 + x_6 + x_7-234), \quad x_1 + x_3 + x_4 + x_6 + x_7 + x_8-234, \\
&117 - x_3 - x_5 - x_8,\quad \tfrac{1}{2}(x_1 + x_3 + 3 x_4 - x_6 - x_7), \quad  2 x_1 - x_3 + x_5,\\
&x_1 + x_3 - x_4 + x_6 - x_7, \quad x_1 + x_3 - x_4 - x_6 + x_7,\quad -2 (x_1 + x_3) + x_8, \\
& 351 - 3 x_1 - x_4 - x_5 - x_6 - x_7 - x_8).
\end{align*}
 
 
Recall that if $E:=\sum_{j\in\{1, 3, 4, 6, 7, 8\}}E_j$, then $\charvec_SE=0$.
Consider $E_j$ where $j\in\{1, 3, 4, 6, 7, 8\}$. Then $E_j=EE_j$ and so $\charvec_SE_j=\charvec_SEE_j=0$.
So we have $(aQ)_j=0$ for $j\in\{1, 3, 4, 6, 7, 8\}$, and hence
\[
2x_4=x_6=x_7 = x_8, \quad x_3 = 117 - x_5 - x_8,\quad x_5 = x_1-\tfrac{3}{2} x_8+117.
\]


Therefore,
\begin{align*}
a&=\left(1,x_1,0,\frac{x_8}{2}-x_1,\frac{x_8}{2},x_1-\frac{3 x_8}{2}+117,x_8,x_8,x_8,-x_1-\frac{5 x_8}{2}+234\right),\\
aQ&=\left(352, 0, 64 ( 2 x_8-117), 0, 0, 32 (117 + 4 x_1 - 2 x_8), 0, 0, 0, 64 (117-2 x_1 - x_8)\right).
\end{align*}

We now take a different fusion scheme yielding a strongly regular graph $\mathcal{F}$.
Let $A= \sum_{i \in \{2,7,8\}}A_i$ where the $A_i$ are adjacency matrices of $\mathcal{A}$, ordered according to the matrix $P$ above. 
Then $A$ is the adjacency matrix of a strongly regular graph $\mathcal{F}$ with parameters $(4096,1638, 662, 650)$.
 The matrix of eigenvalues for $\mathcal{F}$ is
\[
P_\mathcal{F}=\begin{bmatrix}
1&1638&2457\\
1&38&-39\\
1&-26&25
\end{bmatrix}\]
 and the matrix of dual eigenvalues $Q_\mathcal{F}$ is exactly the same as $P_\mathcal{F}$. 
From the inner distribution $a$ for $S$, it follows that 
$|\mathcal{F}(v)\cap S|=2x_8$ for all $v\in S$, where $\mathcal{F}(v)$ denotes the neighbourhood of $v$ in $\mathcal{F}$.
From the outer distribution of $S$, it follows that
$|\mathcal{F}(v)\cap S|=2y_8$ for all $v\notin S$.
Therefore, 
\[
\charvec_SA=2x_8\charvec_S +2y_8 (\ones-\charvec_S)
\]
where $\ones$ is the `all ones' vector,
and so $(2x_8-2y_8-1638)\charvec_S+2y_8\ones$ is an eigenvector for $A$. In particular, $\charvec_S$ is annihilated by one
of the non-principal minimal idempotents of $\mathcal{F}$.
The inner distribution for $S$, with respect to $\mathcal{F}$, is 
$a_\mathcal{F}:=(1,2x_8,351-2x_8)$ and therefore, its MacWilliams transform is
\[
a_\mathcal{F}Q_\mathcal{F}=\left(352, 64 (2x_8-117), 32(351 - 4 x_8)\right).
\]
Since $\charvec_S$ is annihilated by one
of the non-principal minimal idempotents, $(a_\mathcal{F}Q_\mathcal{F})_{j}=0$ for either $j=1$ or $j=2$. 
So there are two cases to consider.

\begin{description}
\item[Case 1: $(a_\mathcal{F}Q_\mathcal{F})_{1}=0$]
Here we have $x_8=117/2$ and so
\begin{align*}
a&=\left(1, x_1, 0, \tfrac{117}{4} - x_1, \tfrac{117}{4}, \tfrac{117}{4} + x_1, \tfrac{117}{2}, \tfrac{117}{2},\tfrac{117}{2}, \tfrac{351}{4} -x_1\right),\\
aQ&=(352, 0, 0, 0, 0, 128 x_1,0,0,0, 32 (117 - 4 x_1)).
\end{align*}
This implies that $S$ is design-orthogonal to the subfield design given in \ref{subfield}. So by Roos' Theorem \ref{roos}, 
\[
|S\cap \F_{4^3}|=\frac{|S||\F_{4^3}|}{|\Gamma|}=\frac{352\cdot 64}{4096}=\frac{11}{2}.
\]
This is a contradiction as $|S\cap \F_{4^3}|$ is an integer.\medskip

\item[Case 2: $(a_\mathcal{F}Q_\mathcal{F})_{2}=0$]
Here we have $x_8=351/4$ and 
\begin{align*}
a&=\left(1, x_1, 0, \tfrac{351}{8} - x_1, \tfrac{351}{8}, x_1-\tfrac{117}{8}, \tfrac{351}{4}, \tfrac{351}{4}, \tfrac{351}{4}, \tfrac{117}{8} - x_1\right),\\
aQ&=\left(352, 0, 3744, 0, 0, 32 (4x_1-\tfrac{117}{2}),0,0,0,32(\tfrac{117}{2}-4x_1)\right).
\end{align*}
In particular, $aQ\ge 0$ implies that $x_1=\tfrac{117}{8}$ and hence
\[
aQ=(352, 0, 3744, 0, 0, 0, 0, 0, 0, 0)
\]
and we have $\charvec_S\in V_0\perp V_2$.
So $S$ is design-orthogonal to $S^\tau$, and so by Roos' Theorem \ref{roos},
\[
|S\cap S^\tau|=\frac{|S||S^\tau|}{|\Gamma|}=\frac{352\cdot 352}{4096}=\frac{121}{4}.
\]
This is a contradiction as $|S\cap S^\tau|$ is an integer.
\end{description}
Both cases lead to a contradiction, and so there is no Delsarte coclique.
\end{proof}

\section{Conclusion}

The existence of a two-character 430-cap would not only yield a new 78-cap of $\PG(5,4)$, but it would yield a two-character cap (see Lemma \ref{twocharacter}),
and hence a new uniformly packed code with parameters $[78,6]$ (over $\F_4$).
So indeed, if it can be shown that the Hill 78-cap of $\PG(5,4)$ is the only two-character cap in this space (up to projectivity), then
a two-character 430-cap does not exist. The following is well-known, but the author cannot find it in print, and so it is proved here.

\begin{lemma}\label{twocharacter}
Let $\mathcal{C}$ be a $430$-cap of $\PG(6,4)$, such that every hyperplane intersects in 78 or 110 elements.
Let $\Sigma$ by a hyperplane intersecting $\mathcal{C}$ in 78 elements. Then $\Sigma\cap \mathcal{C}$
is a cap (of a copy of $\PG(5,4)$) such that every hyperplane of $\Sigma$ intersects in 14 or 22 elements.
\end{lemma}

\begin{proof}
Let $\mu_i$ be the number of elements of $\Sigma\cap \mathcal{C}$ in the $i$-th hyperplane of $\Sigma$. 
Then the usual counting arguments show that  
\begin{align*}
\sum_{i}\mu_i&=78 \cdot  \gaussian{5}{1}_{\mathbb{F}_4}=78\times 341,\\
\sum_{i}\mu_i(\mu_i-1)&=78 \cdot 77\cdot  \gaussian{4}{1}_{\mathbb{F}_4}=78\times 6545,\\
\sum_{i}\mu_i(\mu_i-1)(\mu_i-2)&=78 \cdot 77\cdot 76 \cdot  \gaussian{3}{1}_{\mathbb{F}_4}=78\times 122892.
\end{align*}

Take a hyperplane $H$ of $\Sigma$. Then it is on five hyperplanes, each meeting the $430$-cap $\mathcal{C}$ in at most 110 points.
Since each point of the cap is in at least one of these hyperplanes (n.b., the span of a point and $H$ is a hyperplane of $\PG(6,4)$),
we have
\[
4(110-| H\cap \mathcal{C} |) + 78\ge 430
\]
and hence $| H\cap \mathcal{C} |\le 22$.
So $(\mu_i - 14)^2 (22-\mu_i)\ge 0$ for all $i$.
From the displayed equations above, we have 
\begin{align*}
\sum_i (\mu_i - 14)^2 (22-\mu_i )&= \sum_i \left(-\mu_i(\mu_i-1)(\mu_i-2)+47\mu_i(\mu_i-1) -763\mu_i+4312\right)\\
&= 78\left(-122892+47\times 6545-763\times 341+75460\right)\\
&=0.
\end{align*} 
Therefore, $\mu_i\in\{14,22\}$ as required.
\end{proof}

We remark that the largest coclique of $\Gamma$ that we have been able to find by computation has size 119, but no
well established technique that bounds the size of a coclique seemed to eliminate 352 cocliques immediately.
This includes eigenvalue bounds, spherical code bounds, the No-Homomorphism Lemma, and the Clique-Adjacency polynomial.


\end{document}